\documentclass{article}

\newtheorem{theorem}{Theorem}[section]
\newtheorem{corollary}[theorem]{Corollary}

\newtheorem{lemma}[theorem]{Lemma}
\newtheorem{proposition}[theorem]{Proposition}

\newtheorem{remark}{Remark}[section]

\evensidemargin\oddsidemargin
\begin{document}

\author{Vadim E. Levit and Eugen Mandrescu \\
Department of Computer Science\\
Holon Academic Institute of Technology\\
52 Golomb Str., P.O. Box 305\\
Holon 58102, ISRAEL\\
\{levitv, eugen\_m\}@barley.cteh.ac.il}
\title{Matrices and $\alpha $-Stable Bipartite Graphs}
\date{}
\maketitle

\begin{abstract}
A square $(0,1)$-matrix $X$ of order $n\geq 1$ is called fully
indecomposable if there exists no integer $k$ with $1\leq k\leq n-1$, such
that $X$ has a $k$ by $n-k$ zero submatrix. The reduced adjacency matrix of
a bipartite graph $G=(A,B,E)$ (having $A\cup B=\{a_{1},...,a_{m}\}\cup
\{b_{1},...,b_{n}\}$ as vertex set, and $E$ as edge set), is $%
X=[x_{ij}],1\leq i\leq m,1\leq j\leq n$, where $x_{ij}=1$ if $a_{i}b_{j}\in E
$ and $x_{ij}=0$ otherwise. A stable set of a graph $G$ is a subset of
pairwise nonadjacent vertices. The stability number of $G$, denoted by $%
\alpha (G)$, is the cardinality of a maximum stable set in $G$. A graph is
called $\alpha $-stable if its stability number remains the same upon both
the deletion and the addition of any edge. We show that a connected
bipartite graph has exactly two maximum stable sets that partition its
vertex set if and only if its reduced adjacency matrix is fully
indecomposable. We also describe a decomposition structure of $\alpha $%
-stable bipartite graphs in terms of their reduced adjacency matrices. On
the base of these findings we obtain both new proofs for a number of
well-known theorems on the structure of matrices due to Brualdi (1966),
Marcus and Minc (1963), Dulmage and Mendelsohn (1958), and some
generalizations of these statements. Several new results on $\alpha $-stable
bipartite graphs and their corresponding reduced adjacency matrices are
presented, as well. Two kinds of matrix product are also considered (namely,
Boolean product and Kronecker product), and their corresponding graph
operations. As a consequence, we obtain a strengthening of one Lewin's
theorem claiming that the product of two fully indecomposable matrices is a
fully indecomposable matrix.
\end{abstract}

\section{Introduction}

Throughout this paper $G=(V,E)$ is a simple (i.e., a finite, undirected,
loopless and without multiple edges) graph with vertex set $V=V(G)$ and edge
set $E=E(G)$. If $A$ is a subset of vertices, $G[A]$ is the subgraph of $G$
spanned by $A$, i.e., the graph having $A$ as its vertex set, and containing
all the edges of $G$ connecting vertices of $A$. By $G-W$ we mean either the
subgraph $G[V-W]$ , if $W\subset V(G)$, or the partial subgraph of $G$
obtained by deleting the edges from $W$, whenever $W\subset E(G)$ (we use $%
G-a$, if $W=\{a\}$). If $A,B$ are disjoint subsets of $V$, then $(A,B)$
stands for the set $\{e=ab:a\in A,b\in B,e\in E\}$. The neighborhood of a
vertex $v\in V$, denoted by $N(v)$, is the set of vertices adjacent to $v$.
For any $A\subset V(G)$, we denote $N_{G}(A)=\cup \{N(x):x\in A\}$, or, if
no ambiguity, $N(A)$. A subset $D\subset V(G)$ is said to be $2$-dominating
in $G$ if $\left| N(v)\cap D\right| \geq 2$, for any vertex $v\in V-D$, \cite
{GunHarRall}. A stable set (i.e., a set containing pairwise nonadjacent
vertices) of maximum size will be referred to as a \textit{maximum stable set%
} of $G$. The \textit{stability number} of $G$, denoted by $\alpha (G)$, is
the cardinality of a maximum stable set of $G$. A perfect matching is a set
of non-incident edges of $G$ covering all its vertices.

A bipartite graph is a triple $G=(A,B,E)$, where $E$ is its edge set and $%
\{A,B\}$ is its bipartition; if $|A|=|B|$, then $G$ is called \textit{%
balanced bipartite}. If $A,B$ are as the only two maximum stable sets of $G$%
, then it is a \textit{bistable bipartite graph}.

A graph $G=(V,E)$ is called:

($\mathit{i}$) $\alpha ^{-}$-\textit{stable} if $\alpha (G-e)=\alpha (G)$ is
valid for any $e\in E$, \cite{GunHarRall};

($\mathit{ii}$) $\alpha ^{+}$-\textit{stable} if $\alpha (G+e)=\alpha (G)$
holds for any $e\notin E$, $e=xy$ and $x,y\in V$, \cite{GunHarRall};

($\mathit{iii}$) $\alpha $-\textit{stable} if it is both $\alpha ^{-}$%
-stable and $\alpha ^{+}$-stable, \cite{LevMan}.

Let $G=(A,B,E)$ be a bipartite graph, where $A=\{a_{1},a_{2},...,a_{m}\}$
and also $B=\{b_{1},b_{2},...,b_{n}\}$. Then $G$ can be characterized by its 
\textit{adjacency matrix}, which is a square $(0,1)$-matrix of order $m+n$%
\[
\left[ 
\begin{array}{cc}
O & X \\ 
X^{t} & O
\end{array}
\right] , 
\]
where $X=[x_{ij}],1\leq i\leq m,1\leq j\leq n$, with $x_{ij}=1$ if $%
a_{i}b_{j}\in E$ and $x_{ij}=0$ otherwise. $X$ is called the \textit{reduced
adjacency matrix} of the bipartite graph $G$. Any $(0,1)$-matrix of size $m$
by $n$ is the reduced adjacency matrix of a bipartite graph. If $G$ is
balanced bipartite, then its reduced adjacency matrix is a square $(0,1)$%
-matrix of order $n=|A|=|B|$. The \textit{term rank} $\rho =\rho (X)$ of a $%
(0,1)$-matrix $X$ of size $m$ by $n$ is the maximal number of $1$'s of $X$
with no two of $1$'s on a line (i.e., on a row or on a column). A collection
of $n$ elements of a square $(0,1)$-matrix $X$ of order $n$ is called a 
\textit{diagonal} of $X$ provided no two elements belong to the same row or
column of $X$. A \textit{nonzero diagonal} of $X$ is a diagonal not
containing any $0$'s.

A square $(0,1)$-matrix $X$ of order $n$ is called \textit{partly
decomposable} if $n=1$ and its unique entry is zero, or $n>1$ and there
exists an integer $k$ with $1\leq k\leq n-1$, such that $X$ has a $k$ by $%
n-k $ zero submatrix. A square matrix is \textit{fully indecomposable}
provided it is not partly decomposable, \cite{MarMinc1}. By permuting the
lines of $X$, the partly decomposable matrix $X$ can be written in the form 
\[
\left[ 
\begin{array}{ll}
X_{1} & O \\ 
X_{2} & X_{3}
\end{array}
\right] , 
\]
where $O$ is a zero matrix of size $k$ by $n-k,X_{1}$ and $X_{3}$ are square
matrices of orders $k$ and $n-k$, respectively.

Decomposition structures of $\alpha ^{+}$-stable and $\alpha $-stable
bipartite graphs were first established in Levit and Mandrescu \cite{LevMan1}%
. On the base of these findings we obtain both new proofs for several
well-known theorems on the structure of matrices due to Brualdi \cite
{Brualdi1}, \cite{Brualdi2}, \cite{Brualdi3}, \cite{Brualdi4}, Marcus and
Minc \cite{MarMinc1}, Dulmage and Mendelsohn \cite{DulMend1}, and also some
generalizations of these statements. Some new results on reduced adjacency
matrices of $\alpha $-stable bipartite graphs are presented, as well. For
example, we show that a connected bipartite graph has exactly two maximum
stable sets that partition its vertex set if and only if its reduced
adjacency matrix is fully indecomposable.

The paper is organized as follows: for the sake of self-consistency, section 
$2$ contains a series of results referring to the structure of bistable, $%
\alpha ^{+}$-stable, and $\alpha $-stable bipartite graphs. We use these
findings further, in section $3$, proving some corresponding assertions for
reduced adjacency matrices associated with bipartite graphs. Sections $4$
and $5$ are dealing with two different kinds of matrix product, (namely,
Boolean and Kronecker), and the corresponding graph operations.

\section{$\alpha $-Stable bipartite graphs}

Haynes et al. proved the following theorem, describing stability properties
of general graphs.

\begin{theorem}
\label{th1}\cite{HayLawBriDut} A graph $G$ is:

($\mathit{i}$) $\alpha ^{-}$-stable if and only if each of its maximum
stable sets is a $2$-dominating set in $G$;

($\mathit{ii}$) $\alpha ^{+}$-stable if and only if no pair of vertices is
contained in all its maximum stable sets.
\end{theorem}

Using Theorem \ref{th1}, we proved the following result from \cite{LevMan},
which in particular, is valid for trees, as Gunther et al. show in \cite
{GunHarRall}.

\begin{theorem}
\label{th2}\cite{LevMan} If $G$ is a connected bipartite graph, then the
following assertions are equivalent:

($\mathit{i}$) $G$ is $\alpha ^{+}$-stable;

($\mathit{ii}$) $G$ has a perfect matching;

($\mathit{iii}$) $G$ possesses two maximum stable sets that partition its
vertex set.
\end{theorem}

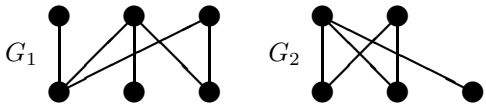
\begin{figure}[h]
\setlength{\unitlength}{1.0cm} 
\begin{picture}(5,1.2)\thicklines

  \multiput(4,0)(1,0){3}{\circle*{0.29}}
  \multiput(4,1)(1,0){3}{\circle*{0.29}}
  \multiput(4,0)(1,0){3}{\line(0,1){1}}
  \put(4,0){\line(1,1){1}}
  \put(4,0){\line(2,1){2}}
  \put(5,1){\line(1,-1){1}}

  \multiput(7.5,0)(1,0){3}{\circle*{0.29}}
  \multiput(7.5,1)(1,0){2}{\circle*{0.29}}
  \multiput(7.5,0)(1,0){2}{\line(0,1){1}}
  \put(7.5,0){\line(1,1){1}}
  \put(7.5,1){\line(1,-1){1}}
  \put(7.5,1){\line(2,-1){2}}

\put(3.5,0.5){\makebox(0,0){$G_{1}$}}
\put(7,0.5){\makebox(0,0){$G_{2}$}}
 \end{picture}
\caption{$\alpha ^{+}$-stable and $\alpha ^{-}$-stable bipartite graphs ${G}%
_{1}$ and $G_{2}$.}
\label{fig1}
\end{figure}

Figure \ref{fig1} illustrates some basic differences between $\alpha ^{+}$%
-stable and $\alpha ^{-}$-stable graphs. Namely, both are bipartite, but $%
G_{1}$ is $\alpha ^{+}$-stable and non-$\alpha ^{-}$-stable (it has a
perfect matching and a non-$2$-dominating maximum stable set), while $G_{2}$
is $\alpha ^{-}$-stable (its unique maximum stable set is $2$-dominating),
and non-$\alpha ^{+}$-stable (it has no perfect matching).

\begin{lemma}
\label{lem1}If $G=(A,B,E)$ is an $\alpha $-stable graph, and $S$ is a
maximum stable set of $G$ meeting both $A$ and $B$, then the subgraph $%
H=G[(S\cap A)\cup (B-S)]$ is $\alpha $-stable.
\end{lemma}

\setlength {\parindent}{0.0cm}\textbf{Proof.} Since $S_{A}=S\cap A$ and $%
B-S_{B}$, (for $S_{B}=S\cap B$), are matched in any perfect matching of $G$, 
$H$ is $\alpha ^{+}$-stable. We show that $H$ is also $\alpha ^{-}$-stable.
Firstly, $S_{A}$ is $2$-dominating, because for any $b\in B-S_{B}$, we have $%
|N(b)\cap S_{A}|=|N(b)\cap S|\geq 2$. $S_{B}$ is also $2$-dominating, since
for any $a\in A-S_{A}$, we have $|N(a)\cap S_{B}|=|N(a)\cap S|\geq 2$. Let $%
X $ be a maximum stable set of $H$, such that both $X_{A}=X\cap A=X\cap
S_{A} $ and $X_{B}=X\cap B=X\cap (B-S_{B})$ are non-empty. $S^{\prime
}=X\cup S_{B} $ is clearly a maximum stable set of $G$, and therefore, we
have: $|N(a)\cap X|=|N(a)\cap X_{B}|=|N(a)\cap S^{\prime }|\geq 2$, for any $%
a\in S_{A}-X_{A}$, and $|N(b)\cap X|=|N(b)\cap X_{A}|=|N(a)\cap S^{\prime
}|\geq 2$, for any $b\in B-S_{B}-X_{B}$, i.e., $X$ is $2$-dominating in $H$.
Consequently, $H$ is $\alpha $-stable, by Theorem \ref{th2}. \rule{2mm}{2mm}%
\setlength
{\parindent}{3.45ex}

\begin{proposition}
\label{prop1}A connected bipartite graph $G$ is $\alpha $-stable if and only
if $G$ can be decomposed as $G=G_{1}\cup G_{2}\cup ...\cup G_{k},k\geq 1$,
such that all $G_{i}=(A_{i},B_{i},E_{i}),1\leq i\leq k$, are
vertex-disjoint, bistable bipartite and $\alpha $-stable.
\end{proposition}

\setlength {\parindent}{0.0cm}\textbf{Proof.} If $G=(A,B,E)$ has $A$ and $B$
as its only two maximum stable sets, then $G$ itself is bistable bipartite
and $\alpha $-stable. Otherwise, let $S$ be a maximum stable set of $G$,
such that both $S_{A}=S\cap A$ and $S_{B}=S\cap B$ are non-empty. By Lemma%
\emph{\ }\ref{lem1}, the subgraphs: $H_{1}=G[(S\cap A)\cup (B-S)]$ and $%
H_{2}=G[(A-S)\cup (S\cap B)]$ are $\alpha $-stable. If they both have only
two maximum stable sets, then they build the decomposition needed.
Otherwise, we continue with this decomposition procedure, until all the
subgraphs we obtain are $\alpha $-stable and have exactly two maximum stable
sets. After a finite number of subpartitions, we get a decomposition $%
G=G_{1}\cup G_{2}\cup ...\cup G_{k},k\geq 1$, such that every $%
G_{i}=(A_{i},B_{i},E_{i}),1\leq i\leq k$, has only $A_{i}$ and $B_{i}$ as
its maximum stable sets.\setlength {\parindent}{3.45ex}

Conversely, let $G=(A,B,E)=G_{1}\cup ...\cup G_{k},k\geq 1$, be such that
each graph $G_{i}=(A_{i},B_{i},E_{i}),1\leq i\leq k$, has only $A_{i}$ and $%
B_{i}$ as its maximum stable sets. Then $G$ is $\alpha ^{+}$-stable, since
it has at least one perfect matching, namely, 
\[
M=\cup \{M_{i}:M_{i}\ is\ a\ perfect\ matching\ in\ G_{i}\}, 
\]
According to\emph{\ }Theorem\emph{\ }\ref{th1}, it suffices to show that any
maximum stable set $S$ of $G$ is also $2$-dominating in $G$. For $S=A$ (and
analogously for $S=B$), suppose $S$ is not $2$-dominating. Hence, there is a
vertex $b\in B_{i}\subseteq B$, such that $|S\cap N(b)|=|\{a\}|=1$. Clearly $%
a\in A_{i}$, and this implies that $A_{i}\cup \{b\}-\{a\}$ is a third
maximum stable set in $G_{i}$, which contradicts the fact that $G_{i}$ is
bistable. Thus, $A_{i}$ (also $B_{i}$) and $S$ are $2$-dominating in $%
G_{i},G $ respectively.

Suppose $S$ meets both $A$ and $B$. We claim that if $i\in \{1,...,k\}$ and $%
S\cap A_{i}\neq \emptyset $, then $A_{i}\subseteq S$ (similarly, if $S\cap
B_{i}\neq \emptyset ,i\in \{1,...,k\}$ then $B_{i}\subseteq S$). Otherwise,
if there is some $j\in \{1,...,k\}$, such that both $S\cap A_{j}$ and $S\cap
B_{j}$ are nonempty, we have:

\[
|S\cap A_{j}|+|S\cap B_{j}|<|A_{j}|=\alpha (G_{j})\ and\ |S\cap
A_{i}|+|S\cap B_{i}|\leq |A_{i}|=\alpha (G_{i}),\ for\ i\neq j. 
\]
Hence, we arrive at the following contradiction:

\[
\alpha (G)=|S|=|S\cap A_{1}|+|S\cap B_{1}|+...+|S\cap A_{k}|+|S\cap
B_{k}|<|A_{1}|+...+|A_{k}|=\alpha (G).
\]
Let $v\in A\cup B-S$. The vertex $v\in B_{i}$ for some $i\in \{1,...,k\}$.
Hence, $S\cap B_{i}=\emptyset $. Consequently, $A_{i}\subseteq S$. and $%
|S\cap N(v)|\geq |A_{i}\cap N(v)|\geq 2$, since $A_{i}$ is $2$-dominating in 
$G_{i}$. Finally, $S$ is also $2$-dominating in $G$, and this completes the
proof. \rule{2mm}{2mm}\newline

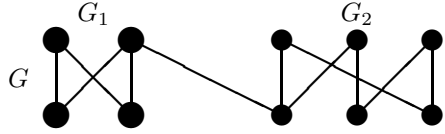
\begin{figure}[h]
\setlength{\unitlength}{1.0cm} 
\begin{picture}(5,1.5)\thicklines

  \multiput(4,0)(1,0){2}{\circle*{0.35}}
  \multiput(4,1)(1,0){2}{\circle*{0.35}}
  \multiput(4,0)(1,0){2}{\line(0,1){1}}
  \put(4,0){\line(1,1){1}}
  \put(4,1){\line(1,-1){1}}
  \put(5,1){\line(2,-1){2}}

  \multiput(7,0)(1,0){3}{\circle*{0.29}}
  \multiput(7,1)(1,0){3}{\circle*{0.29}}

  \multiput(7,0)(1,0){3}{\line(0,1){1}}
  \put(7,0){\line(1,1){1}}
  \put(8,0){\line(1,1){1}}
  \put(7,1){\line(2,-1){2}}

\put(3.5,0.5){\makebox(0,0){$G$}}
\put(4.5,1.35){\makebox(0,0){$G_{1}$}}
\put(8,1.35){\makebox(0,0){$G_{2}$}}
 \end{picture}
\caption{An {example of decomposition: }$G=G_{1}\cup G_{2}$ and $G_{1},G_{2}$
are bistable.}
\label{fig2}
\end{figure}

An example of this decomposition is presented in Figure \ref{fig2}. $%
G=G_{1}\cup G_{2}$ is $\alpha $-stable bipartite and both $G_{1},G_{2}$ are
bistable bipartite.

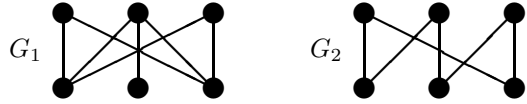
\begin{figure}[h]
\setlength{\unitlength}{1.0cm} 
\begin{picture}(5,1.2)\thicklines

  \multiput(4,0)(1,0){3}{\circle*{0.29}}
  \multiput(4,1)(1,0){3}{\circle*{0.29}}
  \multiput(4,0)(1,0){3}{\line(0,1){1}}
  \put(4,0){\line(1,1){1}}
  \put(4,0){\line(2,1){2}}
  \put(4,1){\line(2,-1){2}}
  \put(5,1){\line(1,-1){1}}

  \multiput(8,0)(1,0){3}{\circle*{0.29}}
  \multiput(8,1)(1,0){3}{\circle*{0.29}}
  \multiput(8,0)(1,0){3}{\line(0,1){1}}
  \put(8,0){\line(1,1){1}}
  \put(9,0){\line(1,1){1}}
  \put(8,1){\line(2,-1){2}}

\put(3.5,0.5){\makebox(0,0){$G_{1}$}}
\put(7.5,0.5){\makebox(0,0){$G_{2}$}}
 \end{picture}
\caption{${G}_{1}$ is non-bistable, ${G}_{2}$ is bistable. }
\label{fig3}
\end{figure}

\begin{theorem}
\label{th3}If $G=(A,B,E)$ is a bipartite graph with at least $4$ vertices,
then the following conditions are equivalent (see examples of a bistable
bipartite graph and a non-bistable bipartite graph in Figure \ref{fig3}):

($\mathit{i}$) $G$ is bistable bipartite;

($\mathit{ii}$) $G$ is $\alpha ^{+}$-stable and $G-a-b$ is $\alpha ^{+}$%
-stable, for any $a\in A$ and $b\in B$;

($\mathit{iii}$) for any $a\in A$ and $b\in B$, $G-a-b$ has a perfect
matching;

($\mathit{iv}$) $G$ is connected and any of its edges is contained in a
perfect matching of $G$;

($\mathit{v}$) $|N(X)|>|X|$ , for any proper subset $X$ of $A$ and of $B$.
\end{theorem}

\setlength {\parindent}{0.0cm}\textbf{Proof.} ($\mathit{i}$) $\Rightarrow $ (%
$\mathit{ii}$) According to Theorem\emph{\ }\ref{th2}, $G$ is $\alpha ^{+}$%
-stable. Let $a\in A,b\in B$ and $H=G-\{a,b\}$. It suffices to show that $%
\alpha (H)=|A-\{a\}|=|B-\{b\}|$.\setlength {\parindent}{3.45ex} Suppose, on
the contrary, that $\alpha (H)=\alpha (G)$; then there is a stable set $S$
in $H$, such that $\alpha (H)=|S|$. Consequently, $S$ is a third maximum
stable set in $G$, in contradiction with the premises on $G$.

($\mathit{ii}$) $\Rightarrow $ ($\mathit{i}$) Clearly, $G$ is connected and $%
\alpha ^{+}$-stable. By Theorem\emph{\ }\ref{th2}, we obtain that $\alpha
(G)=|A|=|B|$. Let $S$ be a third maximum stable set in $G,a\in A-S$ and $%
b\in B-S$. $H=G-\{a,b\}$ is $\alpha ^{+}$-stable and $\alpha
(H)=|A-\{a\}|=|B-\{b\}|=\alpha (G)-1$, by the hypothesis. Since $S$ is
stable in $H$, we obtain the following contradiction $\alpha (G)=|S|\leq
\alpha (H)=\alpha (G)-1$. Consequently, $G$ has only $A$ and $B$ as maximum
stable sets.

($\mathit{ii}$) $\Leftrightarrow $ ($\mathit{iii}$) It is true, according to
Theorem\emph{\ }\ref{th2}.

($\mathit{iii}$) $\Rightarrow $ ($\mathit{iv}$) $G$ is connected, since
otherwise for $a,b$ in different color classes and different connected
components, $G-a-b$ has no perfect matching, contradicting the assumption on 
$G-a-b$. Let $ab$ be an arbitrary edge of $G$ and $M$ be a perfect matching
in $G-a-b$, which exists according to hypothesis. Hence, $M\cup \{ab\}$ is a
perfect matching in $G$ containing $ab$.

($\mathit{iv}$) $\Rightarrow $ ($\mathit{i}$) Suppose, on the contrary, that 
$G$ has a maximum stable set $S$ meeting both $A$ and $B$. If denote $%
S_{A}=S\cap A$ and $S_{B}=S\cap B$, then in any perfect matching of $G$, the
sets $S_{A}$ and $S_{B}$ are matched respectively with $B-S_{B},A-S_{A}$.
Consequently, we obtain that no edge $ab$ joining a vertex $a\in A-S_{A}$
with some vertex $b\in B-S_{B}$ (such an edge must exist, because $G$ is
connected) belongs to some perfect matching of $G$, contradicting the
assumption on $G$. Therefore, $G$ is bistable bipartite.

($\mathit{i}$) $\Rightarrow $ ($\mathit{v}$) Clearly, $\alpha (G)=|A|=|B|$.
Suppose that there is some proper subset $X$ of $A$ such that $|N(X)|\leq
|X| $. Consequently, $(X,B-N(X))=\emptyset $, and hence, $S=X\cup (B-N(X))$
is stable in $G$ with $|S|=|X|+|B-N(X)|\geq |X|+|A-X|=\alpha (G)$. Thus,
since $S$ meets both $A$ and $B$, we infer that $S$ is a third maximum
stable set of $G$, and this is a contradiction, because $G$ is bistable. An
analogous proof can be obtained if $X\subset B$.

($\mathit{v}$) $\Rightarrow $ ($\mathit{i}$) If $|N(X)|>|X|$ holds for any
proper subset $X$ of $A$ and of $B$, it follows that $|A|=|B|\leq \alpha (G)$%
. Assume that some maximum stable set $S$ of $G$ meets both $A$ and $B$.
Then we obtain the following contradiction:

\[
\alpha (G)=|S|=|S\cap A|+|S\cap B|<|N(S\cap A)|+|S\cap B|\leq |B|\leq \alpha
(G). 
\]

Consequently, $G$ must be bistable bipartite. \rule{2mm}{2mm}\newline

The graph $G_{1}$ in Figure \ref{fig3} is non-bistable, since it has $3$
maximum stable sets, but $G_{2}$ is bistable.

\begin{corollary}
\label{cor1}If $G$ is a bistable bipartite graph with at least $4$ vertices,
then 
\[
\cap \{M:M\ is\ a\ perfect\ matching\ of\ G\}=\emptyset .
\]
\end{corollary}

\setlength {\parindent}{0.0cm}\textbf{Proof.} By Theorem\emph{\ }\ref{th3}, $%
G$ is $\alpha ^{+}$-stable, and therefore, it has perfect matchings.
Suppose, on the contrary, that there exists $ab\in \cap \{M:M$ \textit{is a
perfect matching of} $G\}$. If $x\in N(a)-\{b\}$, then, according to Theorem%
\emph{\ }\ref{th3}, $H=G-a-x$ is $\alpha ^{+}$-stable and thus, it has a
perfect matching $M_{0}$, which matches $b$ with some $y\in N(b)-\{a\}$.
Hence, $M_{0}\cup \{ax\}$ is a perfect matching of $G$ and $ab\notin
M_{0}\cup \{ax\}$, contradicting the assumption on $ab$. Therefore, we have $%
\cap \{M:M$ \textit{is a perfect matching of} $G\}=\emptyset $. \rule%
{2mm}{2mm}\setlength {\parindent}{3.45ex}

\begin{proposition}
\label{prop2}A connected bipartite graph $G$ is $\alpha $-stable if and only
if it has perfect matchings and $\cap \{M:M\ is\ a\ perfect\ matching\ of\
G\}=\emptyset $.
\end{proposition}

\setlength {\parindent}{0.0cm}\textbf{Proof.} By Proposition\emph{\ }\ref
{prop1}, $G$ may be decomposed as $G=G_{1}\cup G_{2}\cup ...\cup G_{k},k\geq
1$, such that each $G_{i}=(A_{i},B_{i},E_{i}),1\leq i\leq k$, is bistable
bipartite.\setlength {\parindent}{3.45ex} Taking into account Corollary \ref
{cor1} and the fact that 
\[
\cup \{M_{i}:M_{i}\ is\ a\ perfect\ matching\ of\ G_{i},1\leq i\leq k\} 
\]
is a perfect matching in $G$, we get that $\cap \{M:M\ is\ a\ perfect\
matching\ of\ G\}=\emptyset $.

Conversely, we claim first that from any vertex are issuing at least two
edges contained in some perfect matchings of $G$. Otherwise, there is a
vertex $v$ in $G$, so that only one edge, say $vw$, is contained in a
perfect matching of $G$; such an edge must exist, because $G$ has perfect
matchings. Moreover, since $v$ is matched with a vertex by each such
matching, we infer that $vw$ belongs to all perfect matchings of $G$, in
contradiction with $\cap \{M:M\ is\ a\ perfect\ matching\ of\ G\}=\emptyset $%
. Assume, on the contrary, that $G$ is not $\alpha $-stable, i.e., $G$ is
not $\alpha ^{-}$-stable, since by Theorem\emph{\ }\ref{th2}, $G$ is $\alpha
^{+}$-stable. Therefore, there is a maximum stable set $S$, meeting both $A$
and $B$, and a vertex, say $a\in A$, such that $|N(a)\cap S|=|\{b\}|=1$.
Since from $a$ are issuing at least two edges contained in different perfect
matchings of $G$, we infer that there is at least a vertex $c\in N(a)\cap
(B-S)$, such that $ac$ is in a perfect matching $M$ of $G$. Hence, since $%
|A-S\cap A-\{a\}|<|S\cap B|$, some vertex in $S\cap B$ must be matched by $M$
with some vertex in $S\cap A$, thus contradicting the stability of $S$.
Therefore, $G$ is $\alpha $-stable. \rule{2mm}{2mm}

\begin{proposition}
\label{prop3}A connected balanced bipartite graph $G$ is $\alpha ^{+}$%
-stable if and only if it admits a decomposition as $G=G_{1}\cup ...\cup
G_{k}$, all $G_{i}$ being vertex-disjoint and bistable bipartite.
\end{proposition}

\setlength {\parindent}{0.0cm}\textbf{Proof.} Let $H_{0}=G[M_{0}]$ and $%
H_{1}=G-H_{0}$, where 
\[
M_{0}=\cap \{M:M\ is\ a\ perfect\ matching\ of\ G\}=\emptyset . 
\]

Clearly, $H_{1}$ has $\cap \{M:M\ is\ a\ perfect\ matching\ of\
H_{1}\}=\emptyset $, while $H_{0}$ is either empty or a disjoint union of $%
K_{2}$. According to Propositions\emph{\ }\ref{prop2} and \ref{prop1}, any
connected component of $H_{1}$ has a decomposition in bistable bipartite
subgraphs. Therefore, $G$ admits a decomposition as $G=G_{1}\cup ...\cup
G_{k}$, all $G_{i}$ being vertex-disjoint and bistable bipartite.%
\setlength
{\parindent}{3.45ex}

Conversely, if $G=G_{1}\cup ...\cup G_{k}$, and all $G_{i}$ are bistable
bipartite, then each $G_{i}$ has at least a perfect matching $M_{i}$, and 
\[
\cup \{M_{i}:M_{i}\ is\ a\ perfect\ matching\ of\ G_{i},1\leq i\leq k\} 
\]
is a perfect matching in $G$. Consequently, by Theorem\emph{\ }\ref{th2}, $G$
is $\alpha ^{+}$-stable. \rule{2mm}{2mm}\newline

\begin{figure}[h]
\setlength{\unitlength}{1.0cm} 
\begin{picture}(5,1.5)\thicklines

  \multiput(4,0)(1,0){6}{\circle*{0.29}}
  \multiput(4,1)(1,0){6}{\circle*{0.29}}
  \multiput(4,0)(1,0){6}{\line(0,1){1}}
  \put(4,0){\line(1,1){1}}
  \put(4,1){\line(1,-1){1}}
  \put(5,1){\line(1,-1){1}}

  \multiput(6,0)(1,0){2}{\line(1,1){1}}
  \put(6,1){\line(2,-1){2}}
  \put(8,1){\line(1,-1){1}}
\put(3.5,0.5){\makebox(0,0){$G$}}
\put(4.5,1.35){\makebox(0,0){$G_{1}$}}
\put(7,1.35){\makebox(0,0){$G_{2}$}}
\put(9,1.35){\makebox(0,0){$G_{3}$}}

 \end{picture}
\caption{An example of decomposition into bistable components: $G=G_{1}\cup
G_{2}\cup G_{3}$.}
\label{fig4}
\end{figure}
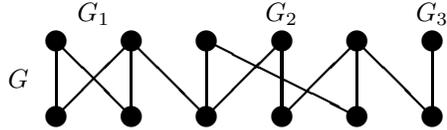

In Figure \ref{fig4} is presented an example of decomposition of an $\alpha
^{+}$-stable bipartite graph into vertex-disjoint and bistable bipartite
components: $G=G_{1}\cup G_{2}\cup G_{3}$.

\section{Matrices and bipartite graphs}

It is not difficult to see that the unity matrix $I_{n},n\geq 1$, is the
reduced adjacency matrix of $nK_{2}$, i.e., of the graph consisting of $n$
disjoint copies of $K_{2}$. Moreover, we have:

\begin{lemma}
\label{lem2}A bipartite graph $G$ is disconnected if and only if its
adjacency matrix $X$ can be written as

\begin{equation}
\emph{\ }\left[ 
\begin{array}{llllll}
X_{1} & O & O & . & . & O \\ 
O & X_{2} & O & . & . & O \\ 
. & . & . & . & . & . \\ 
. & . & . & . & . & . \\ 
. & . & . & . & . & . \\ 
O & O & . & . & . & X_{k}
\end{array}
\right] ,  \label{matrix1}
\end{equation}
where the blocks $X_{1},X_{2},...,X_{k}$ are the adjacency matrices
corresponding respectively to the $k\geq 2$ connected components of $G$.
\end{lemma}

\begin{lemma}
\label{lem3}Let $S$ be a proper subset of the vertex set of graph $G=(A,B,E)$%
, with $p+q$ vertices, where $p=|S\cap A|\geq 1$ and $q=|S\cap B|\geq 1$.
Then $S$ is stable in $G$ if and only if its reduced adjacency matrix $X$
can be written as 
\[
\left[ 
\begin{array}{ll}
X_{1} & O \\ 
X_{2} & X_{3}
\end{array}
\right] ,
\]
where $O$ is a $p$ by $q$ zero matrix.
\end{lemma}

\setlength {\parindent}{0.0cm}\textbf{Proof.} By using an appropriate
indexing for $A$ and for $B$, we may suppose that $S\cap
A=\{a_{1},...,a_{p}\}$ and $S\cap B=\{b_{n-q+1},...,b_{n}\}$. Therefore, $S$
is stable in $G$ if and only if $x_{ij}=0$ for any $i\in \{1,...,p\}$ and $%
j\in \{n-q+1,...,n\}$, i.e., $X$ has exactly the form announced above. \rule%
{2mm}{2mm}\setlength {\parindent}{3.45ex}

\begin{proposition}
\label{prop4}Let $G=(A,B,E)$ be a connected balanced bipartite graph with $2n
$ vertices and $X$ be its reduced adjacency matrix. Then $G$ has a stable
set of $n$ vertices that meets both A and $B$ if and only if $X$ is partly
decomposable.
\end{proposition}

\setlength {\parindent}{0.0cm}\textbf{Proof.} If $p=|S\cap A|$, then $%
q=|S\cap B|=n-p$, and by Lemma\emph{\ }\ref{lem3}, we obtain $X$ in the form 
\[
\left[ 
\begin{array}{ll}
X_{1} & O \\ 
X_{2} & X_{3}
\end{array}
\right] , 
\]
where $O$ is a $p$ by $n-p$ zero matrix, $1\leq p\leq n-1$, i.e., $X$ is
partly decomposable. \rule{2mm}{2mm}\setlength {\parindent}{3.45ex}

\begin{proposition}
\label{prop5}A balanced bipartite graph is bistable if and only if its
reduced adjacency matrix is fully indecomposable.
\end{proposition}

\setlength {\parindent}{0.0cm}\textbf{Proof.} Since a bistable bipartite
graph $G=(A,B,E)$ is connected and has only $A$ and $B$ as maximum stable
sets, Proposition\emph{\ }\ref{prop4} ensures that its reduced adjacency
matrix can not be partly decomposable. The converse is clear. \rule{2mm}{2mm}%
\setlength {\parindent}{3.45ex}\newline

Following the terminology from \cite{DulMend2}, let us recall that for a
balanced bipartite graph $G=(A,B,E)$, a \textit{cover} is a pair of subsets $%
A_{0},B_{0}$ of $A,B$ respectively, such that for every edge $ab\in E$,
either $a\in A_{0}$ or $b\in B_{0}$. $G$ is \textit{cover irreducible} if
its only minimum covers are $\{A,\emptyset \}$ and $\{\emptyset ,B\}$. The
reduced adjacency matrix of a cover irreducible bipartite graph is a \textit{%
cover irreducible matrix}. On the other hand, a bipartite graph $G$ is
called \textit{elementary}, \cite{LovPlum}, if the set, containing any of
its edges that appears in at least one perfect matching, forms a connected
subgraph of $G$. It is shown, \cite{LovPlum}, that elementary bipartite
graphs and the cover irreducible bipartite graphs are the same. It turns out
that bistable bipartite graphs are exactly cover irreducible bipartite
graphs, and fully indecomposable matrices coincide with cover irreducible
matrices. Our approach is based, in principal, on the \textit{bistable
property}. Combining Theorem\emph{\ }\ref{th3} and Proposition\emph{\ }\ref
{prop5}, we get the following result from \cite{Brualdi6}:

\begin{corollary}
\label{cor2}Let $G=(A,B,E)$ be a balanced bipartite graph with $2n$ vertices
and $X$ be its reduced adjacency matrix. Then $X$ is fully indecomposable if
and only if $G$ is connected and any of its edges belongs to a perfect
matching of $G$.
\end{corollary}

We obtain a simple proof for the following characterization of fully
indecomposable matrices, due to Marcus and Minc, \cite{MarMinc1}, and
Brualdi, \cite{Brualdi1}.

\begin{theorem}
\label{th4}A $(0,1)$-matrix $X$ of order $n\geq 2$ is fully indecomposable
if and only if every $1$ of $X$ belongs to a nonzero diagonal and every $0$
of $X$ belongs to a diagonal whose other elements equal $1$.
\end{theorem}

\setlength {\parindent}{0.0cm}\textbf{Proof.} Let $G=(A,B,E)$ be a balanced
bipartite graph with $|A|=|B|=n$, having $X$ as its reduced adjacency
matrix. Then, according to Proposition\emph{\ }\ref{prop5} and\emph{\ }%
Theorem\emph{\ }\ref{th2}, $X$ is fully indecomposable if and only if $G-a-b$
is $\alpha ^{+}$-stable for any $a\in A$ and $b\in B$, i.e., for any $i,j\in
\{1,...,n\}$, the submatrix $Y$, obtained by deleting the row $i$ and the
column $j$ of $X$, has a nonzero diagonal, and this completes the proof. 
\rule{2mm}{2mm}\setlength {\parindent}{3.45ex}\newline

Another consequence is the following result of Marcus and Minc from \cite
{MarMinc1}.

\begin{corollary}
\label{cor3}A fully indecomposable $(0,1)$-matrix $X$ of order $n$ contains
at most $n(n-2)$ zero entries.
\end{corollary}

\setlength {\parindent}{0.0cm}\textbf{Proof.} Let $G=(A,B,E)$ be a balanced
bipartite graph with $X$ as its reduced adjacency matrix. By Proposition%
\emph{\ }\ref{prop5}, $G$ is bistable and according to Theorem\emph{\ }\ref
{th3}.($\mathit{v}$\textit{)}, $|N(v)|\geq 2$ holds for any vertex $v$ of $G$%
. Consequently, any row of $X$ cannot have more than $n-2$ zeros, and hence $%
X$ cannot contain more than $n(n-2)$ zero entries. On the other hand, $%
C_{2n},n\geq 2$, is bistable and its reduced adjacency matrix has exactly $%
n(n-2)$ zero entries. \rule{2mm}{2mm}\setlength {\parindent}{3.45ex}\newline

A $(0,1)$-matrix of order $n\geq 2$ has \textit{total support} provided each
of its $1$'s belongs to a nonzero diagonal. As a consequence, we get the
following result from \cite{Brualdi5}:

\begin{proposition}
\label{prop6}\cite{Brualdi5} Let $X$ be a $(0,1)$-matrix of order $n\geq 2$
with total support, and let $G$ be the bipartite graph whose reduced
adjacency matrix is $X$. Then $G$ is connected if and only if $X$ is fully
indecomposable.
\end{proposition}

\setlength {\parindent}{0.0cm}\textbf{Proof.} Clearly, $X$ is with total
support if and only if any edge of $G$ is contained in a perfect matching of 
$G$. Therefore, taking into account Theorem\emph{\ }\ref{th3} and\emph{\ }%
Proposition \ref{prop5}, we get that: $G$ is connected $\Leftrightarrow $ $G$
is bistable $\Leftrightarrow $ $X$ is fully indecomposable. \rule{2mm}{2mm}%
\setlength {\parindent}{3.45ex}\newline

We can now characterize the bipartite graphs whose reduced adjacency matrix
is with total support.

\begin{proposition}
\label{prop7}The reduced adjacency matrix $X$ of a bipartite graph $G$ has
total support if and only if all connected components of $G$ are bistable
bipartite.
\end{proposition}

\setlength {\parindent}{0.0cm}\textbf{Proof.} If $G$ is connected, then
according to Proposition\emph{\ }\ref{prop5}, $X$ has total support if and
only if $G$ is bistable. If $G$ is disconnected, Lemma\emph{\ } \ref{lem2}
implies that $X$ can be written in the form (\ref{matrix1}), and then $X$
has total support if and only if all the blocks $X_{1},...,X_{k}$ have total
support, i.e., according to Propositions\emph{\ }\ref{prop5} and \ref{prop6}%
, all connected components of $G$ are bistable bipartite. \rule{2mm}{2mm}%
\setlength {\parindent}{3.45ex}

\begin{proposition}
\label{prop8}Let $G$ be a balanced bipartite graph with $2n$ vertices and $X$
be its reduced adjacency matrix. Then the following assertions are
equivalent:

($\mathit{i}$) $G$ is $\alpha ^{+}$-stable;

($\mathit{ii}$) $X$ has a nonzero diagonal;

($\mathit{iii}$) $\rho (X)=n$;

($\mathit{iv}$) $per(X)>0$.
\end{proposition}

\setlength {\parindent}{0.0cm}\textbf{Proof.} By Theorem\emph{\ }\ref{th2}, $%
G$ is $\alpha ^{+}$-stable if and only if it has a perfect matching, i.e.,
its reduced adjacency matrix $X$ has a nonzero diagonal, and this is
equivalent to both ($\mathit{iii}$) and ($\mathit{iv}$). \rule{2mm}{2mm}%
\setlength {\parindent}{3.45ex}\newline

The following result due to Minc is an immediate consequence of the above
proposition.

\begin{corollary}
\label{cor4}\cite{Minc} A $(0,1)$-matrix $X$ of order $n\geq 2$ is fully
indecomposable if and only if every $(n-1)$-square submatrix $Y$ of $X$ has $%
per(Y)>0$.
\end{corollary}

\setlength {\parindent}{0.0cm}\textbf{Proof.} Suppose $X$ is the reduced
adjacency matrix of the balanced bipartite graph $G=(A,B,E)$. According to
Proposition\emph{\ }\ref{prop5}, $X$ is fully indecomposable if and only if $%
G$ is bistable bipartite, and by Theorem\emph{\ }\ref{th3}, this happens if
and only if $G-a-b$ is $\alpha ^{+}$-stable, for any $a\in A$ and $b\in B$,
i.e., by virtue of the Proposition\emph{\ }\ref{prop8}, $per(Y)>0$ holds for
any $(n-1)$-square submatrix $Y$ of $X$. \rule{2mm}{2mm}%
\setlength
{\parindent}{3.45ex}

\begin{theorem}
\label{th5}Let $G$ be a balanced bipartite graph with $2n$ vertices and $X$
be its reduced adjacency matrix. Then $G$ is $\alpha $-stable if and only if 
$X$ can be written as 
\begin{equation}
\emph{\ }\left[ 
\begin{array}{llllll}
X_{1} & X_{12} & X_{13} & . & . & X_{1k} \\ 
O & X_{2} & X_{23} & . & . & X_{2k} \\ 
. & . & . & . & . & . \\ 
. & . & . & . & . & . \\ 
. & . & . & . & . & . \\ 
O & O & O & . & . & X_{k}
\end{array}
\right] ,  \label{matrix2}
\end{equation}
where $X_{1},...,X_{k}$ are fully indecomposable matrices of order at least $%
2$.
\end{theorem}

\setlength {\parindent}{0.0cm}\textbf{Proof.} By Proposition\emph{\ }\ref
{prop1}, $G$ is $\alpha $-stable if and only if it admits a decomposition as 
$G=G_{1}\cup ...\cup G_{k}$, where all $G_{i},1\leq i\leq k$, are
simultaneously $\alpha $-stable and bistable balanced bipartite. Hence,
using an appropriate indexing for the vertices of $G$, $X$ can be written in
the form (\ref{matrix2}), with $X_{1},...,X_{k}$ as reduced adjacency
matrices corresponding to $G_{1},...,G_{k}$, and therefore being fully
indecomposable, by Proposition\emph{\ }\ref{prop5}. Each $X_{i}$ is of order
at least two, since it corresponds to $G_{i}$, which is a bistable bipartite
and $\alpha $-stable graph, i.e., it has at least $4$ vertices. \rule%
{2mm}{2mm}\setlength {\parindent}{3.45ex}

\begin{theorem}
\label{th6}Let $G$ be a balanced bipartite graph with $2n$ vertices and $X$
be its reduced adjacency matrix. Then $G$ is $\alpha ^{+}$-stable if and
only if $X$ can be written in the form (\ref{matrix2}), where $%
X_{1},...,X_{k}$ are fully indecomposable matrices.
\end{theorem}

\setlength {\parindent}{0.0cm}\textbf{Proof.} By Proposition\emph{\ }\ref
{prop1}, $G$ is $\alpha ^{+}$-stable if and only if it admits a
decomposition as $G=G_{1}\cup ...\cup G_{k}$, where all $G_{i},1\leq i\leq k$%
, are bistable balanced bipartite. Hence, using an appropriate indexing for
the vertices of $G$, $X$ can be written in the form (\ref{matrix2}), with $%
X_{1},...,X_{k}$ as reduced adjacency matrices corresponding to $%
G_{1},...,G_{k}$, and therefore being fully indecomposable, by Proposition%
\emph{\ }\ref{prop5}. \rule{2mm}{2mm}\setlength {\parindent}{3.45ex}\newline

As a consequence, we obtain:

\begin{theorem}
\label{th7}(Dulmage and Mendelsohn, \cite{DulMend1}, Brualdi, \cite{Brualdi1}%
). Let $X$ be a $(0,1)$-matrix of order $n$ with term rank $\rho (X)$ equal
to $n$. Then there exist permutation matrices $P$ and $Q$ of order $n$ and
an integer $k\geq 1$ such that $PAQ$ has the form (\ref{matrix2}), where all 
$X_{1},...,X_{k}$ are square fully indecomposable matrices.
\end{theorem}

\setlength {\parindent}{0.0cm}\textbf{Proof.} Let $G$ be a bipartite graph,
whose reduced adjacency matrix is $X$. By Proposition\emph{\ }\ref{prop7}, $%
G $ is $\alpha ^{+}$-stable, and according to Proposition\emph{\ }\ref{prop3}
it admits a decomposition as $G=G_{1}\cup ...\cup G_{k}$, all $G_{i}$ being
bistable bipartite. Hence, using an appropriate indexing for the vertices of 
$G$, $X$ can be written according to the form \textit{(}\ref{matrix2}\textit{%
)}, with $X_{1},...,X_{k}$ as reduced adjacency matrices corresponding to $%
G_{1},...,G_{k}$. Proposition\emph{\ }\ref{prop5} ensures that $%
X_{1},...,X_{k}$ are fully indecomposable. \rule{2mm}{2mm}%
\setlength
{\parindent}{3.45ex}

\begin{corollary}
\label{cor5}Let $X$ be a $(0,1)$-matrix of order $n$ with $\rho (X)=n$. Then
the following assertions are true:

($\mathit{i}$) the intersection of all nonzero diagonals of $X$ is empty if
and only if all $X_{i}$ in the matrix (\ref{matrix2}) are of order at least $%
2$;

($\mathit{ii}$) the number of $1$ by $1$ blocks $X_{i}$ in the matrix (\ref
{matrix2}) is equal to the number of common elements of all nonzero
diagonals of $X$.
\end{corollary}

\begin{corollary}
\label{cor6}(Brualdi, \cite{Brualdi2}) Let $X$ be a square $(0,1)$-matrix of
order $n$ and let $X_{ij}$ denote the matrix obtained from $X$ by striking
the $i$-th row and the $j$-th column. Then $X$ is fully indecomposable if
and only if $per(X_{ij})>0$.
\end{corollary}

\setlength {\parindent}{0.0cm}\textbf{Proof.} Let $G=(A,B,E)$ be a bipartite
graph whose reduced adjacency matrix is $X$. By Proposition\emph{\ }\ref
{prop5}, $X$ is fully indecomposable if and only if $G$ is bistable, i.e.,
for any $a\in A$ and $b\in B,G-a-b$ has a perfect matching (according to
Theorem\emph{\ }\ref{th3}), that is, by Theorem\emph{\ }\ref{th2}, the
matrix $X_{ab}$ has positive permanent. \rule{2mm}{2mm}%
\setlength
{\parindent}{3.45ex}\newline

We ends this section with the following characterization of the reduced
adjacency matrix corresponding to an $\alpha $-stable bipartite graph.

\begin{proposition}
\label{prop9}Let $G$ be a balanced bipartite graph and $X$ be its reduced
adjacency matrix. Then $G$ is $\alpha $-stable if and only if for any
non-zero entry $x_{ij}$ of $X$ there exists a non-zero diagonal of $X$ that
does not contain it.
\end{proposition}

\setlength {\parindent}{0.0cm}\textbf{Proof.} According to Proposition\emph{%
\ }\ref{prop2}, $G$ is $\alpha $-stable if and only if it has perfect
matchings and $\cap \{M:M\ is\ a\ perfect\ matching\ of\ G\}=\emptyset $,
that is $G$ has perfect matchings and for any of its edges $e$ there is a
perfect matching $M$ such that $e\notin M$. In other words, if and only if
for any non-zero entry $x_{ij}$ of $X$ there exists a non-zero diagonal of $%
X $ that does not contain it. \rule{2mm}{2mm}\setlength {\parindent}{3.45ex}

\section{Boolean product of matrices}

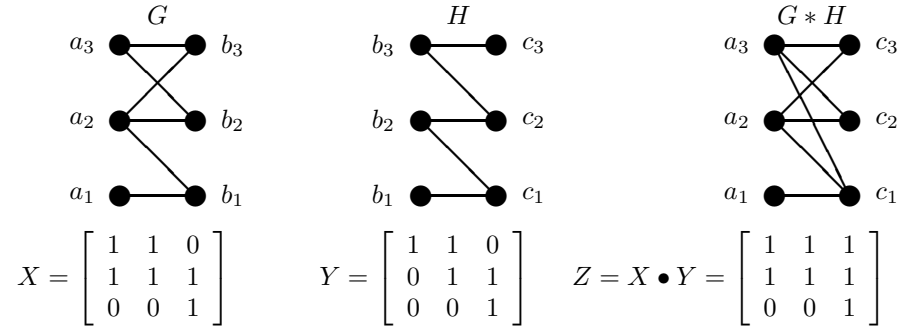
\begin{figure}[h]
\setlength{\unitlength}{1.0cm} 
\begin{picture}(5,4)\thicklines

  \multiput(1.5,1.4)(1,0){2}{\circle*{0.29}}
  \multiput(1.5,2.4)(1,0){2}{\circle*{0.29}}
  \multiput(1.5,3.4)(1,0){2}{\circle*{0.29}}

  \multiput(5.5,1.4)(1,0){2}{\circle*{0.29}}
  \multiput(5.5,2.4)(1,0){2}{\circle*{0.29}}
  \multiput(5.5,3.4)(1,0){2}{\circle*{0.29}}

  \multiput(10.2,1.4)(1,0){2}{\circle*{0.29}}
  \multiput(10.2,2.4)(1,0){2}{\circle*{0.29}}
  \multiput(10.2,3.4)(1,0){2}{\circle*{0.29}}
  
  \multiput(1.5,1.4)(4,0){2}{\line(1,0){1}}
  \multiput(1.5,2.4)(4,0){2}{\line(1,0){1}}
  \multiput(1.5,3.4)(4,0){2}{\line(1,0){1}}

 \multiput(1.5,2.4)(4,0){2}{\line(1,-1){1}}
 \multiput(1.5,3.4)(4,0){2}{\line(1,-1){1}}

  \put(1.5,2.4){\line(1,1){1}}

  \put(10.2,2.4){\line(1,-1){1}}
  \put(10.2,3.4){\line(1,-1){1}}
  \put(10.2,1.4){\line(1,0){1}}
  \put(10.2,2.4){\line(1,0){1}}
  \put(10.2,3.4){\line(1,0){1}}
  \put(10.2,2.4){\line(1,1){1}}
  \put(10.2,3.4){\line(1,-2){1}}

\put(2,3.8){\makebox(0,0){$G$}}
\put(6,3.8){\makebox(0,0){$H$}}
\put(10.7,3.8){\makebox(0,0){$G*H$}}

\put(1,1.4){\makebox(0,0){$a_{1}$}}
\put(1,2.4){\makebox(0,0){$a_{2}$}}
\put(1,3.4){\makebox(0,0){$a_{3}$}}
\put(3,1.4){\makebox(0,0){$b_{1}$}}
\put(3,2.4){\makebox(0,0){$b_{2}$}}
\put(3,3.4){\makebox(0,0){$b_{3}$}}

\put(5,1.4){\makebox(0,0){$b_{1}$}}
\put(5,2.4){\makebox(0,0){$b_{2}$}}
\put(5,3.4){\makebox(0,0){$b_{3}$}}
\put(7,1.4){\makebox(0,0){$c_{1}$}}
\put(7,2.4){\makebox(0,0){$c_{2}$}}
\put(7,3.4){\makebox(0,0){$c_{3}$}}

\put(9.7,1.4){\makebox(0,0){$a_{1}$}}
\put(9.7,2.4){\makebox(0,0){$a_{2}$}}
\put(9.7,3.4){\makebox(0,0){$a_{3}$}}
\put(11.7,1.4){\makebox(0,0){$c_{1}$}}
\put(11.7,2.4){\makebox(0,0){$c_{2}$}}
\put(11.7,3.4){\makebox(0,0){$c_{3}$}}

\put(1.2,0.3){\makebox(0,0){
$\qquad X=\left[ 
\begin{array}{lll}
1 & 1 & 0 \\ 
1 & 1 & 1 \\ 
0 & 0 & 1
\end{array}
\right]  $            }}

\put(5.2,0.3){\makebox(0,0){
$\qquad Y=\left[ 
\begin{array}{lll}
1 & 1 & 0 \\ 
0 & 1 & 1 \\ 
0 & 0 & 1
\end{array}
\right]  $            }}

\put(9.25,0.3){\makebox(0,0){
$\qquad Z=X\bullet Y=\left[ 
\begin{array}{lll}
1 & 1 & 1 \\ 
1 & 1 & 1 \\ 
0 & 0 & 1
\end{array}
\right]  $            }}

 \end{picture}
\caption{The graphs {join operation} and its corresponding Boolean matrix
product.}
\label{fig5}
\end{figure}

Let $G=(A,B,E)$ and $H=(B,C,F)$ be two balanced bipartite graphs on $2n$
vertices. We define the \textit{join} of $G$ with $H$ as the graph $%
P=G*H=(A,C,W)$, where $ac\in W$ if and only if there is $b\in B$, such that $%
ab\in E$ and $bc\in F$. The \textit{Boolean matrix product} of two $(0,1)$%
-matrices $X,Y$ is a $(0,1)$-matrix denoted by $X\bullet Y$ and having the
same zero and non-zero entries as the usual matrix product $XY$; the term 
\textit{Boolean} refers actually to the property of the Boolean addition
operation: $1+1=1$; for an example, see Figure \ref{fig5}). Using this
notation we have the following:

\begin{lemma}
\label{lem4}If $X,Y$ are respectively, the reduced adjacency matrices of the
balanced bipartite graphs $G$ and $H$, then the Boolean matrix product $%
Z=X\bullet Y$ is the reduced adjacency matrix of the graph $P=G*H$.
\end{lemma}

\setlength {\parindent}{0.0cm}\textbf{Proof.} If $X=(x_{ij}),Y=(y_{ij})$ and 
$Z=(z_{ij})$, then clearly we have:

$z_{ij}=\sum_{k=1}^{n}x_{ik}y_{kj}$ $\neq 0\Leftrightarrow $ there exists $%
k\in \{1,...,n\}$ such that $x_{ik}=y_{kj}=1$ $\Leftrightarrow $ there is
some $b_{k}\in B$, so that $a_{i}b_{k}\in E$ and $b_{k}c_{j}\in F$ $%
\Leftrightarrow $ $a_{i}c_{j}\in W$. \rule{2mm}{2mm}%
\setlength
{\parindent}{3.45ex}

\begin{remark}
$X\bullet Y$ is fully indecomposable if and only if $XY$ is fully
indecomposable.
\end{remark}

\begin{corollary}
\label{cor7}Any balanced bipartite graph $G$ on $2n,n\geq 1$, vertices is
isomorphic to $G*nK_{2}$.
\end{corollary}

\begin{proposition}
\label{prop10}Let $G=(A,B,E)$ and $H=(B,C,F)$ be balanced bipartite graphs.

($\mathit{i}$) If $G$ and $H$ are $\alpha ^{+}$-stable, then $G*H$ is $%
\alpha ^{+}$-stable.

($\mathit{ii}$) If one of $G,H$ is $\alpha ^{+}$-stable and the other is
bistable bipartite, then $G*H$ is bistable bipartite.

($\mathit{iii}$) If $G$ and $H$ are bistable bipartite, then $G*H$ is also
bistable bipartite.
\end{proposition}

\setlength {\parindent}{0.0cm}\textbf{Proof. }($\mathit{i}$) Taking into
account the definition of $*$-operation, it is clear that $G*H$ has a
perfect matching, whenever both $G$ and $H$ have a perfect matching. Hence,
Theorem \ref{th2} implies that $G*H$ is $\alpha ^{+}$-stable whenever $G$
and $H$ are both $\alpha ^{+}$-stable.\setlength {\parindent}{3.45ex}

($\mathit{ii}$) Suppose that $G$ is $\alpha ^{+}$-stable and $H$ is bistable
bipartite. If $D$ is an arbitrary proper subset of $A$ or of $C$, then
according to Theorem\emph{\ }\ref{th3} and Hall's marriage theorem we get: $%
|D|<|N_{G}(D)|\leq |N_{H}(N_{G}(D))|=|N_{G*H}(D)|$, i.e., $G*H$ is bistable,
by virtue of the same Theorem\emph{\ }\ref{th3}.

The assertion ($\mathit{iii}$) is a consequence of ($\mathit{ii}$). \rule%
{2mm}{2mm}

\begin{corollary}
\label{cor8}Let $X,Y$ be $(0,1)$-matrices of order $n$. If $per(X)>0$ and $Y$
is fully indecomposable, then $XY$ is fully indecomposable.
\end{corollary}

\begin{corollary}
\label{cor9}(Lewin, \cite{Lewin}) The product of any finite number of fully
indecomposable matrices is a fully indecomposable matrix.
\end{corollary}

\setlength {\parindent}{0.0cm}\textbf{Proof.} Clearly, it is sufficient to
prove the statement for two matrices, say $X$ and $Y$. Let $G=(A,B,E)$ and $%
H=(B,C,F)$ be balanced bipartite graphs, having $X,Y$ respectively, as
reduced adjacency matrices. Lemma\emph{\ }\ref{lem4} implies that $X\bullet
Y $ is the reduced adjacency matrix of the graph $G*H$. By Proposition\emph{%
\ }\ref{prop5}, $G$ and $H$ are bistable bipartite, and according to
Proposition\emph{\ }\ref{prop10}, $G*H$ is also bistable bipartite. Hence,
Proposition\emph{\ }\ref{prop5} ensures that $X\bullet Y$ is fully
indecomposable. Therefore, $XY$ is fully indecomposable, as well. \rule%
{2mm}{2mm}\setlength
{\parindent}{3.45ex}

\begin{corollary}
\label{cor10}(Marcus and Minc, \cite{MarMinc1}) If $X$ is a fully
indecomposable $(0,1)$-matrix, then $XX^{t}$ is fully indecomposable.
\end{corollary}

\section{Kronecker product of matrices}

Let $G=(A,B,E)$ and $H=(C,D,F)$ be two balanced bipartite graphs on $2n$
vertices. The \textit{Kronecker product} of graphs $G$ and $H$ is the graph $%
K=G\otimes H=(A\times C,B\times D,U)$, where $(a,c)(b,d)\in U$ if and only
if $ab\in E$ and $cd\in F$. In these notations we have the following:

\begin{lemma}
\label{lem5}If $X,Y$ are respectively, the reduced adjacency matrices of the
balanced bipartite graphs $G$ and $H$, then the Kronecker matrix product $%
Z=X\otimes Y$ is the reduced adjacency matrix of the graph $K=G\otimes H$.
\end{lemma}

\setlength {\parindent}{0.0cm}\textbf{Proof.} If $X=(x_{ij}),Y=(y_{ij})$ and 
$Z=(z_{ij})$, then we have:

$z_{ij}=z_{(k-1)m+p,(r-1)m+q}=x_{kr}y_{pq}=1$ $\Leftrightarrow $ $x_{kr}=1$
and $y_{pq}=1$

$\Leftrightarrow $ $a_{k}b_{r}\in E$ and $c_{p}d_{q}\in F$ $\Leftrightarrow $
$(a_{k},c_{p})(b_{r},d_{q})\in U$, i.e., $Z$ is the reduced adjacency matrix
of $K$. \rule{2mm}{2mm}\setlength {\parindent}{3.45ex}

\begin{proposition}
\label{prop11}If $G=(A,B,E)$ and $H=(C,D,F)$ are $\alpha ^{+}$-stable, then
their Kronecker product $K=G\otimes H$ is also $\alpha ^{+}$-stable.
\end{proposition}

\setlength {\parindent}{0.0cm}\textbf{Proof.} Let $\{(a_{i},b_{i}):1\leq
i\leq n\}$ and $\{(c_{j},b_{j}):1\leq j\leq m\}$ be perfect matchings in $%
G,H $ respectively, which exist by virtue of Theorem\emph{\ }\ref{th2}.
Hence, according to the same theorem, $K$ is also $\alpha ^{+}$-stable,
since $\{(a_{i},c_{j})(b_{i},d_{j}):1\leq i\leq n,1\leq j\leq m\}$ is a
perfect matching of $K$. \rule{2mm}{2mm}\setlength
{\parindent}{3.45ex}

\begin{corollary}
\label{cor11}Let $X,Y$ be two $(0,1)$-matrices of order $n,m$, respectively.
Then 
\[
\rho (X\otimes Y)\geq \rho (X)\rho (Y),\ and\ if\ \rho (X)=n,\rho (Y)=m,\
then\ \rho (X\otimes Y)=\rho (X)\rho (Y).
\]
\end{corollary}

\setlength {\parindent}{0.0cm}\textbf{Proof.} Let $G=(A,B,E)$ and $H=(C,D,F)$
be bipartite graphs having $X,Y$ as reduced adjacency matrices,
respectively. If the edge sets 
\[
\{a_{i}b_{i}:1\leq i\leq \rho (X)\}\ and\ \{c_{j}b_{j}:1\leq j\leq \rho
(Y)\} 
\]
are maximum matchings in $G,H$ respectively, then 
\[
M=\{(a_{i},c_{j})(b_{i},d_{j}):1\leq i\leq \rho (X),1\leq j\leq \rho (Y)\} 
\]
is a matching in $G\otimes H$, and consequently $\rho (X\otimes Y)\geq
|M|\geq \rho (X)\rho (Y)$. If $\rho (X)=n$ and $\rho (Y)=m$, i.e., both $G$
and $H$ have perfect matchings, then $M$ is a perfect matching in $G\otimes
H $, and this ensures that $\rho (X\otimes Y)=\rho (X)\rho (Y)$. \rule%
{2mm}{2mm}\setlength {\parindent}{3.45ex}

\begin{proposition}
\label{prop12}If $G=(A,B,E)$ is $\alpha $-stable and $H=(C,D,F)$ is $\alpha
^{+}$-stable, then their Kronecker product $K=G\otimes H$ is $\alpha $%
-stable.
\end{proposition}

\setlength {\parindent}{0.0cm}\textbf{Proof.} Let $X,Y,Z$ be the
corresponding reduced adjacency matrices of $G,H$ and $K$. By Proposition%
\emph{\ }\ref{prop9}, for any non-zero entry $%
z_{ij}=z_{(k-1)m+p,(r-1)m+q}=x_{kr}y_{pq}$ of $Z$, there is a non-zero
diagonal $\{x_{1i_{1}},x_{2i_{2}},...,x_{ni_{n}}\}$ of $X$ that does not
contain $x_{kr}$, and clearly the blocks $%
\{x_{1i_{1}}Y,x_{2i_{2}}Y,...,x_{ni_{n}}Y\}$ contain one non-zero diagonal
of $Z$, since $Y$ has at least a non-zero diagonal. According to Proposition%
\emph{\ }\ref{prop9}, $K$ is $\alpha $-stable. \rule{2mm}{2mm}%
\setlength
{\parindent}{3.45ex}

\begin{corollary}
\label{cor12}The Kronecker product of two $\alpha $-stable bipartite graphs
is $\alpha $-stable.
\end{corollary}

In \cite{Brualdi4}, Brualdi proves that:

\begin{theorem}
\label{th8}The Kronecker product of two fully indecomposable matrices is a
fully indecomposable matrix.
\end{theorem}

As a consequence, we get:

\begin{corollary}
The Kronecker product of two bistable bipartite graphs is a bistable
bipartite graph.
\end{corollary}

\section{Conclusions}

In this paper we investigated the intimate relationship existing between the
structure of both $\alpha ^{+}$-stable and $\alpha $-stable bipartite
graphs, and the structure of their corresponding reduced matrices. The
mutual transfer of the results was done via the following bridge: \textit{%
bistable\ bipartite\ graphs\ vis-a-vis\ fully\ indecomposable\ matrices}.

On the base of this duality, we have obtained new proofs and extensions of
several well-known theorems on matrices, and on the other hand, new
characterizations of $\alpha ^{+}$-stable or $\alpha $-stable bipartite
graphs.

\end{document}